\documentclass[letterpaper, 10 pt, conference]{ieeeconf}  

\IEEEoverridecommandlockouts                              

\usepackage{graphics} 
\usepackage{epsfig} 
\usepackage{mathptmx} 
\usepackage{times} 
\usepackage{amsmath} 
\usepackage{amssymb}  
\usepackage{cleveref}
\usepackage{algorithm}
\usepackage{algpseudocode}
\usepackage{xcolor}
\usepackage{csquotes}

\newcommand{\vect}[1]{\boldsymbol{#1}}
\newcommand{\mylist}[1]{\overline{\MakeUppercase{#1}}}

\title{\LARGE \bf
Bi-Level Route Optimization and Path Planning\\with Hazard Exploration
}

\author{Jimin Choi$^{1}$, Grant Stagg$^{2}$, Cameron K. Peterson$^{2}$, Max Z. Li$^{3}$
\thanks{*This work was supported by NSF IUCRC Phase I: Center for Autonomous Air Mobility and Sensing (CAAMS) Award No. 2137195.}
\thanks{$^{1}$Jimin Choi is with the Department of Aerospace Engineering, University of Michigan,
        Ann Arbor, MI 48109, USA; 
        {\tt\small jiminch@umich.edu}}%
\thanks{$^{2}$Grant Stagg and Cameron K. Peterson are with the Department of Electrical Engineering, Brigham Young University,
        Provo, UT 84602, USA; 
        {\tt\small \{ggs24, cammy.peterson\}@byu.edu}}%
        \thanks{$^{3}$Max Z. Li is with the Department of Aerospace Engineering, Department of Civil and Environmental Engineering, and Department of Industrial and Operations Engineering, University of Michigan, Ann Arbor, MI 48109, USA; 
        {\tt\small maxzli@umich.edu}}%
}

\begin{document}

\maketitle
\thispagestyle{empty}
\pagestyle{empty}

\begin{abstract}
Effective risk monitoring in dynamic environments such as disaster zones requires an adaptive exploration strategy to detect hidden threats. We propose a bi-level unmanned aerial vehicle (UAV) monitoring strategy that efficiently integrates high-level route optimization with low-level path planning for known and unknown hazards. At the high level, we formulate the route optimization as a vehicle routing problem (VRP) to determine the optimal sequence for visiting known hazard locations. To strategically incorporate exploration efficiency, we introduce an edge-based centroidal Voronoi tessellation (CVT), which refines baseline routes using pseudo-nodes and allocates path budgets based on the UAV's battery capacity using a line segment Voronoi diagram. At the low level, path planning maximizes information gain within the allocated path budget by generating kinematically feasible B-spline trajectories. Bayesian inference is applied to dynamically update hazard probabilities, enabling the UAVs to prioritize unexplored regions. Simulation results demonstrate that edge-based CVT improves spatial coverage and route uniformity compared to the node-based method. Additionally, our optimized path planning consistently outperforms baselines in hazard discovery rates across a diverse set of scenarios.
\end{abstract}

\section{Introduction}

Effective hazard detection by autonomous agents is crucial, particularly in environments such as disaster zones and contaminated areas, where human responders face significant risks or access constraints. Unmanned aerial vehicles (UAVs) provide a viable alternative, enabling real-time monitoring and assessment in such scenarios~\cite{MOHDDAUD202230}. As a result, there has been increasing interest in deploying UAVs for critical missions, including disaster response, industrial inspection, and environmental hazard detection~\cite{a2ghazard, MISHRA20201}. Recent studies have explored various UAV-assisted monitoring strategies, including routing optimization, cooperative UAV networks, and data-driven approaches~\cite{ ERDELJ201772, wang_reinforcement_2019, choi2025bandit}.

However, previous research predominantly focuses on static environments with predefined hazards, making conventional vehicle routing problems (VRPs) inadequate for dynamically emerging threats~\cite{CALAMONERI2024123766}. Real-world disasters and hazardous areas are inherently dynamic, with new hazards arising unpredictably---examples of such dynamic, new hazards include sudden structural collapses, gas leaks, and secondary fires. In such scenarios, UAVs must balance exploration and exploitation, efficiently monitoring known hazards while actively searching for unforeseen dangers. While most existing methods rely on static models to optimize UAV routes, they often overlook the need for real-time adaptability, resulting in limited responsiveness to newly emerging threats. 
This study addresses this gap by incorporating adaptive exploration mechanisms, enabling UAVs to dynamically search for unknown hazards during path planning while maintaining efficient monitoring of known hazards.

Furthermore, a critical limitation of current UAV-based monitoring approaches is the separate treatment of routing and path planning rather than a joint approach. 
While routing determines the sequence in which hazard locations are visited, path planning ensures safe and feasible navigation between them. Most existing studies have addressed routing and path planning as independent problems, and only a limited number of works have attempted to integrate them~\cite{Kiesel_Burns_Wilt_Ruml_2012, 9216819}. However, previous efforts to jointly consider routing and path planning remain constrained by their focus on static environments with predefined observation points, and lack of an exploration perspective. 
To address this gap, we propose an integrated framework that optimizes UAV routing for effective hazard monitoring while incorporating explorative path planning to generate feasible UAV trajectories. This approach ensures efficient and responsive hazard detection in complex, evolving environments.

\section{Contributions of Work}
The contributions of this research are threefold. First, to overcome the limitations of treating routing and path planning as separate problems, we propose a novel bi-level framework that integrates high-level routing optimization with low-level path planning. This unified approach captures the interdependence of the two tasks while satisfying both operational and dynamic constraints. Second, building on this framework, we introduce a routing strategy that supports efficient exploration in environments with emerging threats. Instead of limiting UAVs to fixed observation points, the strategy diversifies route coverage and efficiently allocates the remaining flight resources to support exploration. Third, we develop a risk-aware path planning method that accounts for potential unknown hazards through Bayesian inference, while generating kinematically feasible B-spline trajectories to complement the routing strategy.

\section{Problem Statement}

\begin{figure}[htbp]
	\begin{center}
		\includegraphics[width= 0.65\linewidth]{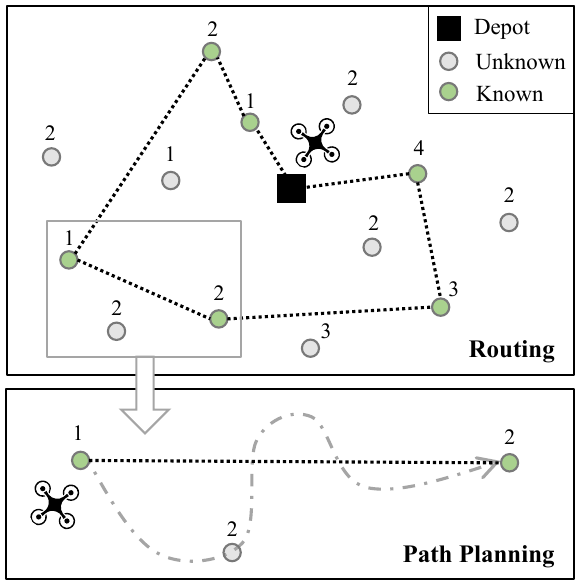}
	\end{center}
	\vspace{-0.5cm}
	\caption{Concept of operations for route optimization in the presence of known hazards and exploring unknown hazards using path planning. Green circles indicate known hazards, gray circles represent unknown hazards, and the depot is the UAV’s base for charging, data transfer, or coordination.}
	\label{fig:conops}
\end{figure}

We propose a bi-level routing and path planning framework for UAV-based monitoring missions that accounts for both known and unknown hazards. \Cref{fig:conops} illustrates the concept of operations. First, we optimize the routing phase to determine the optimal sequence for UAVs to visit known hazards. However, as a graph-based combinatorial optimization problem, VRP defines only a high-level route between nodes and does not consider UAV dynamics. Thus, a dedicated path planning phase refines the routes generated at the routing level by incorporating UAV dynamics and exploration objectives, ensuring that the planned paths remain consistent with routing strategies rather than being designed in isolation. In this path-planning phase, UAVs follow the high-level route while inserting kinematically feasible paths. We define any remaining battery capacity---interpreted as additional travel distance---as the marginal budget, which is used to explore potential unknown hazards. This approach maximizes hazard detection while ensuring compliance with operational constraints.

\begin{figure}[htbp]
  \centering
  \includegraphics[scale=0.45]{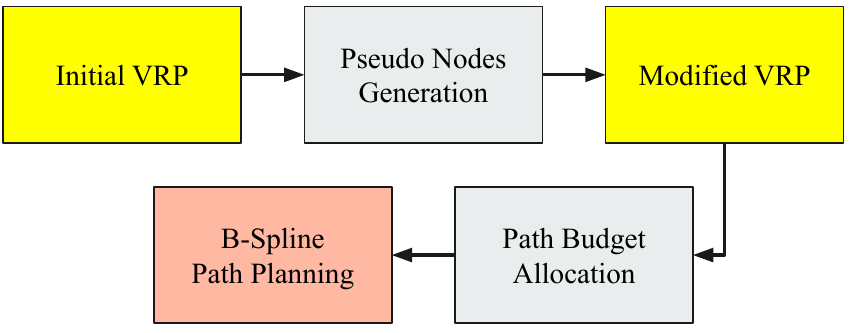}
  \caption{Bi-level routing and path planning process: Yellow blocks represent route optimization, determining the optimal sequence for visiting hazards. Gray blocks adjust the route to improve exploration efficiency in the path planning phase. The red block represents path planning, which generates feasible UAV trajectories considering dynamics.}
  \label{fig:bilevel_flow}
\end{figure}

\Cref{fig:bilevel_flow} illustrates the proposed methodology. The routing phase first determines the optimal sequence for visiting known hazards. To improve coverage and support future exploration, we introduce pseudo-nodes—--virtual waypoints added between known hazards to guide UAVs through less-visited areas. These are placed based on the overall structure of the route to fill in coverage gaps. After route optimization, the remaining path budget is distributed across the route segments. In the path planning phase, UAVs use this path budget to generate kinematically feasible trajectories that follow the planned route and actively explore regions likely to contain unknown hazards.


\subsection{High Level: Vehicle Routing Problem with Known Hazards}

VRP is a well-known optimization problem focused on determining efficient routes while considering constraints. We tackle a UAV-based VRP for hazard monitoring, ensuring that all hazard nodes are visited. Operational constraints like battery life, payload, and time windows are reflected in the path budget and can be adjusted through the problem formulation. Since battery endurance is the primary constraint in UAV missions~\cite{uav_battery}, we model the VRP with distance constraints, using battery capacity as the path budget.

Let \( N_k \) be the number of known hazard nodes and define \( I_k = \{ 1, 2, \dots, N_k \} \) as their corresponding index set. The index \( 0 \) represents the depot node, which serves as the UAV's physical base for activities such as battery recharging, data transfer, and mission coordination. The set of all nodes, including the depot, is given by \( I_0 = I_k \cup \{ 0 \} \). The distance between any two nodes \( i, j \in I_0 \) is denoted as \( d_{ij} > 0 \). 

The set of UAVs is denoted as \( M \), where each UAV \( m \in M \) has a maximum allowable travel distance \( D_m \), determined by its battery capacity. A binary variable \( y_{ijm} \in \{0,1\} \) is used to indicate whether UAV \( m \) travels directly from node \( i \) to node \( j \). To eliminate subtours, we introduce a continuous auxiliary variable \( u_{im} \in \mathbb{R}_+ \), representing the relative order in which UAV \( m \) visits node \( i \in I_k \). 
The mathematical formulation of our VRP for high-level route optimization used for both the initial and modified VRPs is:
\begin{subequations} \label{eq:vrp}
\begin{align}
    \min \quad & \sum_{m \in M} \sum_{i \in I_0} \sum_{j \in I_0} d_{ij} y_{ijm} \label{eq:vrp_obj} \\
    \text{subject to:} \quad & \nonumber \\
    & \sum_{m \in M} \sum_{\substack{j \in I_0 \\ j \neq i}} y_{ijm} = 1, \quad \forall i \in I_k, \label{eq:vrp_nodevisit} \\
    & \sum_{\substack{j \in I_0 \\ j \neq i}} y_{ijm} = \sum_{\substack{j \in I_0 \\ j \neq i}} y_{jim}, \quad \forall m \in M, \forall i \in I_0, \label{eq:vrp_flow} \\
    & \sum_{i \in I_0} \sum_{j \in I_0} d_{ij} y_{ijm} \leq D_m, \quad \forall m \in M, \label{eq:vrp_battery} \\
    & \sum_{j \in I_k} y_{0 j m} = 1, \quad \forall m \in M, \label{eq:vrp_start} \\
    & \sum_{i \in I_k} y_{i 0 m} = 1, \quad \forall m \in M, \label{eq:vrp_end} \\
    & u_{0m} = 0, \quad \forall m \in M, \label{eq:vrp_mtz_base} \\
    & 1 \le u_{im} \le |I_k|, \quad \forall i \in I_k,\; \forall m \in M, \label{eq:vrp_mtz_bounds} \\
    & u_{jm} \ge u_{im} + 1 - |I_k|(1 - y_{ijm}), \nonumber \\
    &\forall m \in M,\; \forall i \neq j \in I_k, \label{eq:vrp_mtz_constraint}\\
    & y_{ijm} \in \{0,1\}, \quad \forall i,j \in I_0, \forall m \in M, \label{eq:vrp_integer}\\
    & u_{im} \in \mathbb{R}_+, \quad \forall i \in I_k, \forall m \in M. \label{eq:vrp_u_domain}
\end{align}
\end{subequations}

In the formulation, \Cref{eq:vrp_obj} defines the objective function, minimizing the total travel distance (or battery consumption) by summing \( d_{ij} \) for each segment where \( y_{ijm} = 1 \). \Cref{eq:vrp_nodevisit} ensures each known hazard node \( i \in I_k\) is visited exactly once, guaranteeing full coverage. \Cref{eq:vrp_flow} enforces flow conservation at every node \( i \in I_0 \), ensuring each UAV \( m \) departs as many times as it arrives, maintaining route continuity. \Cref{eq:vrp_battery} limits each UAV’s total travel distance to its battery capacity \( D_m \), preventing energy overuse. \Cref{eq:vrp_start,eq:vrp_end} require each UAV to start and end its route at the depot (node \( {0} \)). \Cref{eq:vrp_mtz_base,eq:vrp_mtz_bounds,eq:vrp_mtz_constraint} implement subtour elimination using the Miller–Tucker–Zemlin (MTZ) formulation, which enforces a consistent visit order and prevents disconnected cycles. Finally, \Cref{eq:vrp_integer,eq:vrp_u_domain} define decision variables and domains, including the binary routing variable \( y_{ijm} \) and the continuous auxiliary variable \( u_{im} \) for subtour control.

\subsection{Low Level: Path Planning for Unknown Hazards Exploration}~\label{sec:path_planning}
The goal of the path planner is to minimize the posterior probability of an unknown hazard existing in the region. This is achieved using an objective function that seeks to minimize the average posterior probability evaluated on a grid of points in the region. Each section of the path---defined between nodes from the VRP---is planned separately.

We define the set of known hazard locations as \( \mylist{x}_k = \{\vect{x}_{k,1},\hdots,\vect{x}_{k, N_k}\} \), with \(N_k\) denoting the number of known hazards. Similarly, the set of unknown hazards is given by $\mylist{x}_u = \{\vect{x}_{u,1},\hdots,\vect{x}_{u,N_u}\}$, with $N_u$ representing the number of unknown hazards. 
We assume that the locations of the unknown hazards are spatially correlated with the locations of the known hazards.  Specifically, the probability of an unknown hazard occurring at location $\vect{x}_u$ given a single known hazard at $\vect{x}_k$ is given by
\begin{align}
    P(\vect{x}_u \mid\vect{x}_k) 
    &= e^{-\lambda \|\vect{x}_u - \vect{x}_k\|^2_2},
\end{align}
where $\lambda$ is the decay rate for how fast the spatial correlation drops with distance. For this work, we model the spatial correlation with this Gaussian kernel function; however, other differentiable spatial correlation functions could be employed depending on the specific application. Additionally, we assume that there is an underlying prior probability of $P_h$ for an unknown hazard existing at any location. 

Combining the underlying prior with the probability of an unknown hazard existing at $\vect{x}_u$ given all known hazards $\mylist{X}_k$, we get a prior unknown hazard probability of 
\begin{equation} \label{eq:combined_prior}
    P(\vect{x}_u \mid  \mylist{x}_k) = 1 - \left(    \prod_{i=1}^{N_k} (1 -P(\vect{x}_u \mid \vect{x}_{k,i}))(1-P_h)\right).
\end{equation}
The prior given in \Cref{eq:combined_prior} assumes that the locations of the unknown hazards are spatially correlated with the known hazards. If there is no spatial correlation, or if the correlation structure is unknown, a uniform prior of $0.5$ could be used. 

\begin{figure}[t]
  \centering
\includegraphics[scale=0.6]{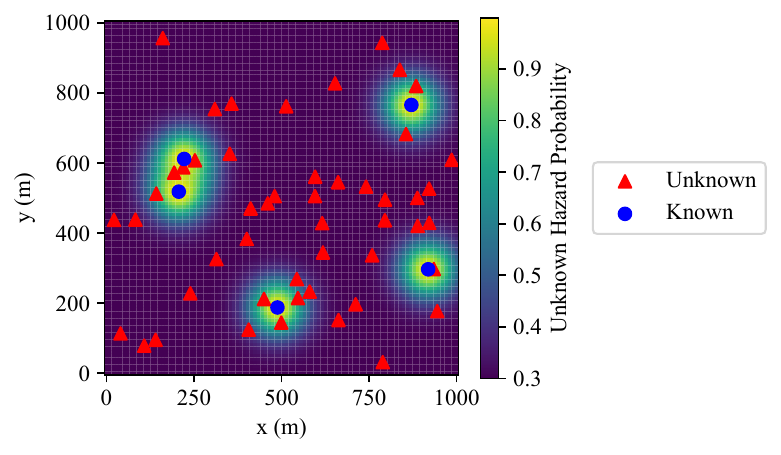}
  \caption{This figure shows the unknown hazard probability from five known hazards and a baseline prior of $ P=0.3$. The known hazards are shown in blue, and the unknown hazards are in red. One hundred unknown hazards are randomly selected using this distribution.}
  \label{fig:hazard_distribution}
\end{figure}

\Cref{fig:hazard_distribution} shows an example of the unknown hazard probability.
Using a Bayesian approach, we will find the probability of an undiscovered hazard existing at a location given a planned path in Section~\ref{sec:path_planning}. We will then update the path to minimize that posterior probability evaluated across points in the region to find a trajectory between nodes. 

As the agent\footnote{We use \enquote{agent} and \enquote{UAV} interchangeably throughout this paper.} travels through the environment, it takes measurements at a frequency of $f_s = 1/\Delta_s$. We assume that it successfully measures the unknown hazard at location $\vect{x}_u$ according to the following probability:
\begin{equation}\label{eq:measurement_model}
P(\vect{x}_{s} \mid \vect{x}_{u}) = e^{-\beta \|\vect{x}_{s}-\vect{x}_u\|^2/_2},
\end{equation}
where $\beta$ is a steepness parameter that determines the effective sensing range of the agent.
We also assume there is some probability the agent detects a hazard when none is present (false alarm), given as
\begin{equation}
    P(\vect{x}_{s} \mid \neg\vect{x}_{u}) = P_{f_a}.
\end{equation}

These probabilistic models characterize the agent’s sensing capabilities in the environment. To ensure that planned trajectories are physically realizable, we now introduce a kinematic model that constrains the agent's motion. We assume the agent follows a unicycle kinematic model
\begin{equation} \label{eq:dynamics}
\begin{bmatrix}
\dot{x}(t) \\
\dot{y}(t) \\
\dot{\theta}(t)
\end{bmatrix} = 
\begin{bmatrix}
v(t)\cos{\theta(t)} \\
v(t)\sin{\theta(t)} \\
u(t)
\end{bmatrix},
\end{equation}
where $x(t)$ is the agent's position east of the origin, $y(t)$ is the agent's position north of the origin, and $\theta(t)$ is the agent's heading measured counter-clockwise from the east axis. The agent can control its velocity $v(t)$ and its turn-rate $u(t)$

To generate trajectories that satisfy the unicycle dynamics, we require a path representation that is both smooth and easily differentiable. One effective choice is B-splines, which provide an efficient and flexible parameterization of continuous paths.
B-splines are defined by a set of control points \( \mylist{C} = (\vect{c}_1, \vect{c}_2, \dots, \vect{c}_{N_c}) \), where \( N_c \) is the number of control points, and a sequence of knot points \( \vect{t}_k = (t_0 - k\Delta_t, \dots, t_0 - \Delta_t, t_0, t_0 + \Delta_t, \dots, t_f, t_f + \Delta_t, \dots, t_f + k\Delta_t) \), where \( t_0 \) and \( t_f \) denote the start and end times of the trajectory, respectively; \( k \) is the order of the B-spline, \( \Delta_t = (t_f - t_0)/N_k \) defines the time spacing between knot points, and \( N_k \) is the number of internal knot points. The trajectory $\vect{p}(t)$ is then expressed as a weighted sum of basis functions 
\begin{equation} \label{eq:b-spline_basis}
\vect{p}(t) = \sum_{i=1}^{N_c} B_{i,k}(t) \vect{c}_i,
\end{equation}
where \( B_{i,k}(t) \) are basis functions derived using the Cox-de Boor recursive formula~\cite{cox1972numerical}. B-splines are particularly advantageous for path planning due to the finite support of their basis functions, meaning that modifications to control points at one end of the spline do not impact the other. Additionally, they are straightforward to differentiate, making them useful for applications requiring smooth trajectories.

\section{Routing Strategy}
Our proposed routing strategy considers exploration in path planning, rather than relying solely on known hazard nodes. It includes two key features during route optimization: augmenting the route with pseudo-nodes to improve spatial coverage, and allocating the remaining path budget to support efficient local exploration.

\subsection{Pseudo-nodes Generation for Route Augmentation}
Exploration efficiency is a key factor in the route optimization phase of path planning. To enhance UAV exploration and improve path coverage, we introduce pseudo-nodes. Initially, the optimized route considers all of the known nodes, establishing a visit sequence. However, if known nodes are unevenly distributed, large inter-node gaps may result in inefficient coverage by an agent traversing only to those known nodes. In practice, since the known nodes are randomly placed, they often cluster in specific regions, leaving other areas underexplored or entirely uncovered. This leads the VRP solution---which seeks the shortest route satisfying constraints between known nodes---to be biased toward these clustered known nodes, reducing overall exploration potential by neglecting sparsely covered areas that may contain hazards. To mitigate this, pseudo-nodes are incorporated into the route optimization process in a manner (described below) that provides more uniformly distributed waypoints to the path planner.  This, in turn, guides the UAV through less-visited areas and ensuring more effective overall coverage.

Pseudo-nodes are generated based on the spatial distribution of known nodes using centroidal Voronoi tessellation (CVT), a specialized Voronoi diagram wherein the generating points coincide with the center of mass of their respective Voronoi cells~\cite{CVT_1999}. In a CVT, a given set of generators \( x_i \) partitions the domain \( \Omega \) into Voronoi cells \( V_i \), where each cell consists of all points closest to its respective generator
\begin{equation}
    V_i = \{ x \in \Omega \mid \| x - x_i \|_2 \leq \| x - x_j \|_2, \forall j \neq i \}.
\end{equation}
To achieve a CVT, each generator \( x_i \) must be the centroid  of its Voronoi cell, computed as:
\begin{equation}
    x_i = \frac{\int_{V_i} x \rho(x) d x}{\int_{V_i} \rho(x) d x},
\end{equation}
where \( \rho(x) \) is a density function weighing the importance of different regions. A uniform \( \rho(x) \) results in evenly distributed Voronoi regions, which we call node-based CVT in this paper.

We introduce two key modifications and refer to this modified version as edge-based CVT. First, the centroid calculation is modified by introducing a non-uniform density function that depends on the distance to edges in the VRP route, thereby reshaping the Voronoi cells. Specifically, we define the density function as
\begin{equation}
    \rho(x) = \alpha + \beta \, d(x, \mylist{E}),
\end{equation}
where \( d(x, \mylist{E}) \) represents the Euclidean distance from the point \( x \) to the nearest edge in the set \( \mylist{E}=\{e_1,e_2,...,e_m\}\), where each edge \( e_k \) connects two nodes \( i, j \in I_0 \) in the VRP route. The scalar parameters \(\alpha, \beta >0\) control the sensitivity of the density with respect to edge proximity. This modification biases the centroid computation by assigning greater weight to points that are farther from edges. As a result, the computed centroid is pulled away from edge-dense areas, reshaping the Voronoi cells accordingly. 

The second modification is made in the generator update process. In our approach, generators are divided into two types: fixed generators, corresponding to known hazard nodes, and adaptive generators, corresponding to pseudo-nodes.  Pseudo-nodes are placed at the center of edge-based Voronoi regions that are not associated with any known hazard node. Each known hazard remains fixed and serves as the centroid of its own Voronoi cell, while only the pseudo-nodes are updated through the CVT process. To compute centroids under the modified density function, a set of sample points \( \mylist{s} \subset \Omega \) is drawn from the domain. Each sample is assigned to its nearest generator to define the Voronoi regions used in the CVT update.
Algorithm~\ref{alg:edge_based_CVT} formalizes the edge-based CVT process using the modified density function and a fixed or adaptive generator structure.


\begin{algorithm}[htbp]
\caption{Edge-Based CVT pseudo-nodes Generation}
\label{alg:edge_based_CVT}
\begin{algorithmic}[1]
    \Require Known nodes $\mylist{x}_k$, edges $\mylist{e}$, domain $\Omega$, Number of pseudo-nodes $N_{\text{pseudo}}$, weight parameters $\alpha, \beta$
    \State Initialize pseudo-nodes $\mylist{x}_p$ (e.g., randomly within $\Omega$)
    \For{iteration = $1$ to max\_iter}
        \State Sample points $\mylist{s}$ in $\Omega$ and assign to nearest generator
        \For{each pseudo-node $x_i \in \mylist{X}_p$}
            \State Compute the weighted centroid:
            \[
            x_i = \frac{\sum_{s \in \mylist{s}_i} s \rho(s)}{\sum_{s \in \mylist{s}_i} \rho(s)}, \quad \rho(s) = \alpha + \beta d(s, \mylist{E})
            \]
        \EndFor
    \EndFor
    \State Re-solve VRP with $\mylist{x}_k \cup \mylist{x}_p$
    \If{any UAV route exceeds battery capacity}
        \State Adjust or remove pseudo-nodes and re-solve VRP
    \EndIf
    \State \textbf{Return:} Combined set of known nodes and pseudo-nodes: $\mylist{x}_k \cup \mylist{x}_p$
\end{algorithmic}
\end{algorithm}

\subsection{Marginal Path Budget Allocation for Exploration} \label{sec:initial_budget_allocation}
Once the pseudo-nodes are generated, we compute an optimized route that visits both known nodes and pseudo-nodes. This route allows the UAV to visit all known nodes while efficiently exploring for unknown hazards. Using our VRP formulation, the optimized route does not consume the entire budget, and some margin is left to ensure feasibility. However, since the path budget is allocated to the entire route rather than to individual edges, managing the remaining budget effectively becomes challenging. 
To address this, we distribute the marginal budget across individual VRP edges, allowing the agent to use the leftover budget for exploration. 

Let \( \mylist{x}_a = \mylist{x}_k \cup \mylist{x}_p = \{\vect{x}_{a,i}, \vect{x}_{a,2}, \dots, \vect{x}_{a, N_a}\} \) denote the augmented node set, consisting of both known hazard nodes \( \mylist{x}_k \) and pseudo-nodes \( \mylist{x}_p \). Each \( \vect{x}_{a,i} \in \mathbb{R}^2 \) represents the spatial position of a node.  A line segment Voronoi diagram is a generalization of the classic Voronoi diagram where the generators are line segments, rather than points~\cite{line_segment_voronoi}. We employ this diagram to partition the area and proportionally allocate the path budget to each edge before passing the finalized allocation to the path planner. 
The diagram assigns each point in \( \Omega \) to the closest edge from a given set of edges \( \mylist{E}_{final} \), which is obtained from the updated VRP solution over the augmented node set \( \mylist{x}_a \). Each edge \( e_k \in \mylist{E}_{\text{final}} \) connects two nodes \( \vect{x}_{a,i}, \vect{x}_{a,j} \in \mylist{x}_a \) and is represented as \( e_k = (\vect{x}_{a,i}, \vect{x}_{a,j})\). The corresponding Voronoi region for edge $e_k$ is given by:
\begin{equation}
V_k = \{ x \in \Omega \mid d(x, e_k) \leq d(x, e_l), \forall l \neq k \},
\end{equation}
where the distance function \( d(x, e_k) \) is computed as
\begin{equation}
d(x, e_k) =
\begin{cases}
\| x - \vect{x}_{a,i} \|_2, & \text{if } t < 0, \\
\| x - \vect{x}_{a,j} \|_2, & \text{if } t > 1, \\
\| x - (\vect{x}_{a,i} + t (\vect{x}_{a,j} - \vect{x}_{a,i})) \|_2, & \text{otherwise}.
\end{cases}
\end{equation}
Note that \( t \) is the projection factor given by
\begin{equation}
t = \frac{(x - \vect{x}_{a,i}) (\vect{x}_{a,j} - \vect{x}_{a,i})}{\| \vect{x}_{a,j} - \vect{x}_{a,i} \|_2^2}.
\end{equation}

In the context of path planning, the line segment Voronoi diagram distributes the initial given path budget \( B \) across the edges proportionally to their assigned Voronoi regions. Given the area \( A_k \) corresponding to each edge \( e_k \), the budget allocation is computed as
\begin{equation}
B_k = B \frac{A_k}{\sum_{j} A_j},
\end{equation}
where \( A_k \) is the total area of the grid cells assigned to edge \( e_k \). By using area-based allocation, edges assigned larger Voronoi regions receive a proportionally greater share of the path budget, supporting more efficient exploration.

\section{Path Planning Strategy}
The path planner aims to determine a route between the nodes assigned by the VRP that maximizes information gain while adhering to the limited path budget $B_k$ allocated to that segment in Section~\ref{sec:initial_budget_allocation}. Rather than generating the entire route in a single step, the planner optimizes each inter-node segment individually.
To accomplish this, we present a novel path planning algorithm that generates kinematically feasible trajectories while leveraging prior information from known hazards for efficient exploration. Given the locations of known hazards and planned sampling points, we estimate the probability that an unknown hazard exists in that region. We then adjust the path to minimize the average probability of an unknown hazard existing in the region, considering both known hazard locations and the agent's traversal path.

Using the prior distribution of an unknown hazard probability (Equation~\eqref{eq:combined_prior}) and the measurement model (Equation~\eqref{eq:measurement_model}), we wish to find the probability of an unknown hazard being present at a location $\vect{x}_u$ after the agent traverses a proposed path $\vect{p}(t)$. We note that the probability of not detecting an unknown hazard is $P(\neg\vect{x}_{s} \mid  \vect{x}_u ) = 1-P(\vect{x}_{s} \mid \vect{x}_u).$ 
The probability of not detecting a hazard at a location given that there is no hazard is $P(\neg\vect{x}_s \mid \neg\vect{x}_u) = 1-P_{f_a}$.
We assume that detection events at different sampling locations are conditionally independent given the unknown hazard location $\vect{x}_u$. Let $\mylist{X}_s = \{\vect{x}_{s,1}, \hdots, \vect{x}_{s,N_s}\}$ denote the set of sampling locations, where $N_s$ is the number of locations.
Then, the probability of \textit{no detection} at any of the sampling locations, given a hazard at $\vect{x}_u$, is  
\begin{equation}\label{eq:combined_no_detection_given_hazard}
    P\left(\bigcap_{j=1}^{N_s} \neg D_j \mid \vect{x}_u \right)= \prod_{j=1}^{N_s} \left(1 - P(D_j \mid \vect{x}_u)\right),
\end{equation}
where $D_j$ is the event of detecting the hazard at location $\vect{x}_{s,j}$.

Similarly, the combined probability of receiving false detections at the set of locations $\mylist{x}_s$ is
\begin{equation}
    P\left(\bigcap_{\vect{x}_{s,j}\in\mylist{x}_s} \neg\vect{x}_{s,j} \mid \neg \vect{x}_u \right)= (1 - P_{f_a})^{N_s}.
\end{equation}

Using this information, the posterior probability that an unknown hazard exists at $\vect{x}_u$ given that the agent sampled at $\mylist{x}_s$ and did not detect that hazard is
\begin{align}\label{eq:posterior_prob}
    \Phi\left(\vect{x}_u,\mylist{X}_s,\mylist{x}_k \right)
    &=
    P\left(\vect{x}_u \mid  \bigcap_{\vect{x}_{s,j}\in\mylist{x}_s} \neg\vect{x}_{s,j}\right) \notag \\ &=   \frac{P\left(\bigcap_{\vect{x}_{s,j}\in\mylist{x}_s} \neg\vect{x}_{s,j} \mid \vect{x}_u \right)P(\vect{x}_u)}{P\left(\bigcap_{\vect{x}_{s,j}\in\mylist{x}_s} \neg\vect{x}_{s,j}\right)},
\end{align}
where $P(\bigcap_{\vect{x}_{s,j}\in\mylist{x}_s} \neg\vect{x}_{s,j} \mid\vect{x}_u)$ is found in Equation~\eqref{eq:combined_no_detection_given_hazard}, $P(\vect{x}_u)$ is the prior probability of an unknown hazard and can be found using Equation~\eqref{eq:combined_prior} if unknown hazards are spatially correlated with known hazards, and 
\begin{align}
    P\left(\bigcap_{\vect{x}_{s,j}\in\mylist{x}_s} \neg\vect{x}_{s,j}\right)
    = \; &P\left(\bigcap_{\vect{x}_{s,j}\in\mylist{x}_s} \neg\vect{x}_{s,j} \mid \vect{x}_u\right)P(\vect{x}_u) \notag \\
    + \; &P\left(\bigcap_{\vect{x}_{s,j}\in\mylist{x}_s} \neg\vect{x}_{s,j} \mid \neg \vect{x}_u\right)(1 - P( \vect{x}_u)).
\end{align}

We evaluate the posterior probability (Equation~\eqref{eq:posterior_prob}) at a discretized set of points evenly distributed across the Voronoi cell $V_k$ corresponding to the current edge $e_k$ being planned.
These grid points are denoted as $\mylist{x}_{g,k}$ for the $k^\textit{th}$ edge, and there are $N_g$ grid points. 
The objective function (Equation~\eqref{eq:path_objective}) is defined as the average posterior probability evaluated at  
$\mylist{x}_{g,k}$. We assume the agent collects samples along its trajectory $\vect{p}(t)$ at discrete times $\vect{t}_s = (t_0, t_0 + \Delta_s, \hdots, t_f)$, where $\Delta_s = (t_f - t_0)/N_s$ is the time between samples. Let $\mylist{x}_s$ denote the set of past sampling locations from previous edges.
Given a set of known hazard locations $\mylist{x}_k$, the objective function is  
\begin{equation}\label{eq:path_objective}
    \Gamma(\vect{p}(t),\mylist{x}_g,\mylist{x}_s,\mylist{x}_k) = \frac{1}{N_g}\sum_{\vect{x}_{g,i} \in \mylist{x}_g} \Phi\big(\vect{x}_{g,i},\mylist{X}_s \cap \vect{p}(\vect{t}_s), \mylist{x}_k\big),
\end{equation}
where $\Phi(\vect{x}_{g,i}, \cdot)$ represents the posterior probability at grid point $\vect{x}_{g,i}$, conditioned on the current sampling locations along the trajectory $\vect{p}(\vect{t}_s)$, the past sampling locations $\mylist{x}_s$, and the known hazard locations $\mylist{x}_k$.


We seek to plan each segment of the trajectory (between nodes from the VRP) to minimize the posterior probability The path planning problem is the following constrained optimization problem: 
\begin{subequations} \label{eq:optimization}
\begin{align}
    \mylist{C}_{k}, t_{f_{k}} &= \operatorname*{argmin}_{\mylist{C}_k,t_{f_{k}}} \; \Gamma(\vect{p}(t),\mylist{x}_g,\mylist{x}_s,\mylist{x}_k) \label{eq:optimization_util} \\
    \text{subject to} \quad 
    & \vect{p}(0) = \vect{x}_0,\; \vect{p}(t_f) = \vect{x}_f, \label{eq:pos_constraints} \\
    & \vect{\dot{p}}(0) = \vect{v}_0,\; \vect{\dot{p}}(t_f) = \vect{v}_f, \label{eq:vel_constraints} \\
    & v_{lb} \leq v(\vect{t}_s) \leq v_{ub}, \label{eq:velocity_constraint} \\
    & u_{lb} \leq u(\vect{t}_s) \leq u_{ub}, \label{eq:turn_rate_constraint} \\
    & -\kappa_{ub} \leq \kappa(\vect{t}_s) \leq \kappa_{ub}, \label{eq:curve_constraint} \\
    & L(\vect{p}(t)) = B_k. \label{eq:path_length_constraint}
\end{align}
\end{subequations}

The constraints of the path planner are broken down as follows. \Cref{eq:pos_constraints} ensure that the path starts and ends at the nodes of the edge that is being planned. \Cref{eq:vel_constraints} constrains the path to have a specific starting and ending velocity. The starting velocity is set to the current velocity of the agent, and the final velocity is set so that the heading of the agent is pointed towards the following node from the VRP. The next three constraints (\Cref{eq:velocity_constraint,eq:turn_rate_constraint,eq:curve_constraint}) ensure that the path is kinematically feasible. 
The velocity of the trajectory can be found as $v(t) = ||\dot{\vect{p}}(t)||_2.$ The turn rate is $u(t) = (\dot{\vect{p}}(t) \times \ddot{\vect{p}}(t))/||\dot{\vect{p}}(t)||_2^2.$ The final kinematic constraint is the curvature, and can be found as $\kappa(t) = u(t)/v(t).$ The final constraint, \Cref{eq:path_length_constraint}, is the path budget constraint. The function $L(\vect{p}(t))$ returns the path length of the B-spline. This is computed through numerical integration.  

\section{Simulation and Results}
In this section, we present simulation results validating the proposed bi-level framework. First, we compare node-based and edge-based CVT strategies for pseudo-nodes generation in high-level route optimization, evaluating their effect on spatial coverage. Then, we assess the impact of different path planning strategies on the number of unknown hazards discovered during exploration.

\subsection{High-Level Route Optimization}
\begin{figure*}[htbp]
  \centering
  \includegraphics[scale=0.6]{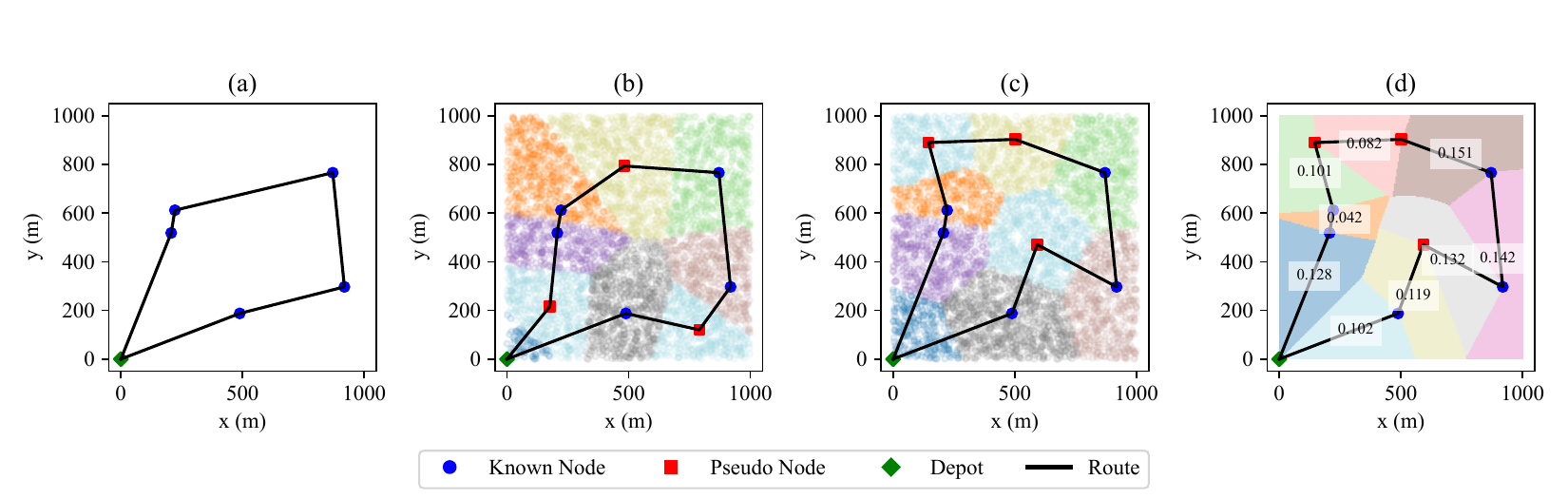}
  \caption{Example of a routing strategy simulation: (a) Initial vehicle routing solution using only known nodes. (b) Routing solution incorporating additional pseudo-nodes using node-based CVT. (c) Routing solution incorporating additional pseudo-nodes using edge-based CVT. (d) Initial path budget allocation visualized through a line segment Voronoi diagram.}
  \label{fig:routing_strategy_sim}
\end{figure*}
\begin{figure*}[htbp]
  \centering
  \includegraphics[scale=0.6, trim=0 0 0 30, clip]{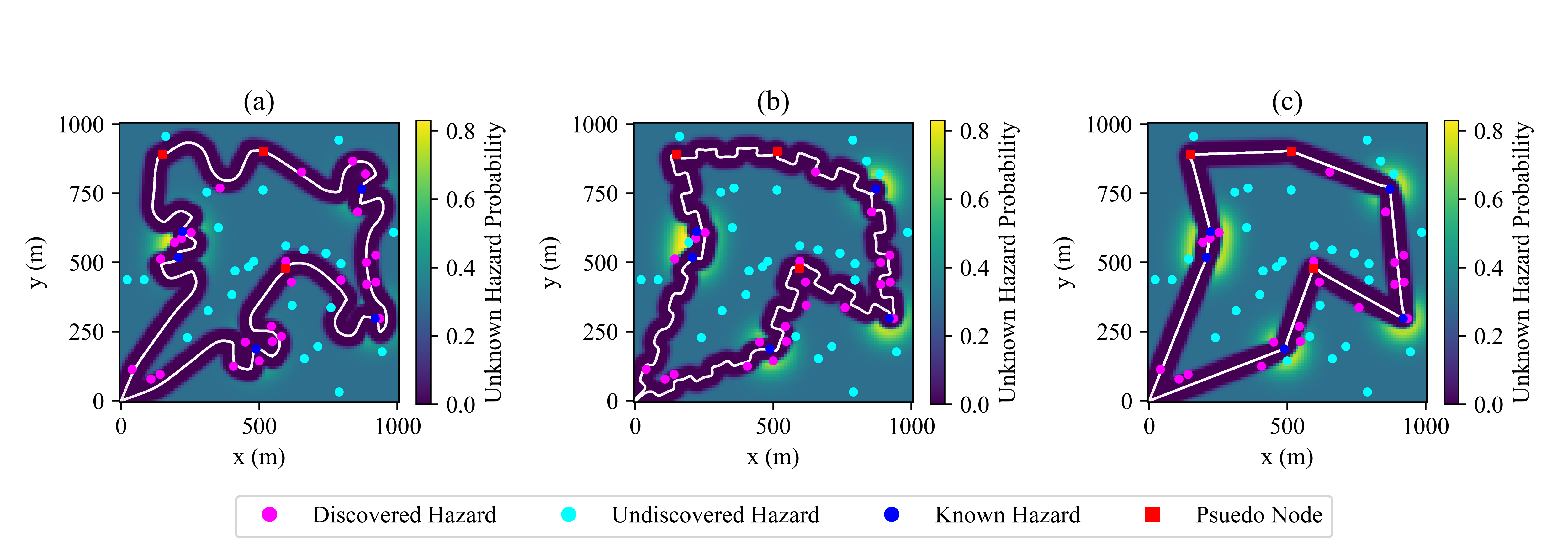}
  \caption{A comparison of the three path planning methods: (a) shows our optimization approach; (b) is the ``lawnmower" pattern; and (c) is straight-line paths between nodes. The posterior probability of discovering hazards is shown, as well as the discovered (magenta) and undiscovered (cyan) hazards. The white line indicates the actual path flown by UAV.}
  \label{fig:planned_paths}
\end{figure*}

We conducted simulations to compare the efficiency of node-based and edge-based CVT in pseudo-nodes placements. The experiments were performed in a bounded 2D domain \(\Omega = [0,1000] \times [0,1000] \;m^2\) with randomly generated known nodes under a total budget constraint of \(B_{total}=5000 \;m\). VRP was solved using the PyVRP library, which is based on a genetic algorithm, to establish a baseline route~\cite{pyvrp}. To compare performance, we defined two metrics. The first metric, edge coverage ratio (ECR), measures the proportion of the domain covered by edges:
\begin{equation}
    \text{ECR} = N_c / N_g,
\end{equation}
where \( N_c \) is the number of grid cells intersected by at least one edge, and \( N_g \) is the total number of grid cells. The domain was discretized into a \(50 \times 50\) uniform grid (\( N_g = 2500 \)), and an edge was considered to intersect a cell if any segment of the VRP path passed through it.  

The second performance metric is the edge density variance (EDV), which quantifies the uniformity of edge distribution:
\begin{equation}
    \text{EDV} = \frac{1}{N_g} \sum_{i=1}^{N_g} (X_i - \bar{X})^2, \quad \bar{X} = \frac{1}{N_g} \sum_{i=1}^{N_g} X_i,
\end{equation}
where \( X_i \) is the number of edges passing through grid cell \( i \), and \( \bar{X} \) is the mean number of edges per cell. A lower EDV indicates a more uniform edge distribution, which is preferred for balanced exploration. It suggests that edges are spread more evenly across the domain rather than being concentrated in specific regions.

Simulations were repeated 100 times with randomly varying node distributions, adjusting the number of known and pseudo-nodes to be placed. The number of known nodes ranged between 5 and 15, while the number of pseudo-nodes varied between 1 and 5. For edge-based CVT, we set \(\alpha = 0.1\) and \(\beta = 0.9\) to prioritize pseudo-nodes in areas with fewer edges. 

\Cref{table_cvt_comparison} shows that edge-based CVT consistently achieves better spatial coverage. This demonstrates that edge-aware CVT enhances pseudo-node placement efficiency, leading to a more uniform route distribution while ensuring the feasibility of VRP.

\Cref{fig:routing_strategy_sim} illustrates an example of the routing strategy. Subfigure~\ref{fig:routing_strategy_sim}(a) presents the known nodes and the optimal route connecting them. The UAV departs from the depot, visits all known nodes, and returns to the depot. Subfigures~\ref{fig:routing_strategy_sim}(b) and~\ref{fig:routing_strategy_sim}(c) compare the results of the node-based and edge-based CVT. Both methods improve coverage compared to the initial VRP solution, but the edge-based CVT is more effective as it accounts for UAV movement along edges. The different colors in the regions represent the final clustered areas obtained through the CVT processes. After incorporating pseudo-nodes generated by edge-based CVT, VRP is solved again, and subfigure~\ref{fig:routing_strategy_sim}(c) illustrates the resultant modified routing outcome. Compared to subfigure~\ref{fig:routing_strategy_sim}(a), which only considers known nodes, the updated solution provides improved coverage, facilitating a more favorable setting for path planning and exploration. Subfigure~\ref{fig:routing_strategy_sim}(d) shows the initial path budget allocation based on the modified routing, where the numbers in each Voronoi cell represent a normalized path budget summing to $1$. The partitioned regions further guide the path planner, highlighting potential areas for exploration along given edges.

\begin{table}[htbp]
\caption{Performance Comparison of CVT Methods}
\label{table_cvt_comparison}
\begin{center}
\begin{tabular}{|c||c|c|}
\hline
\textbf{Method} & \textbf{ECR} & \textbf{EDV} \\
\hline
Original Known Nodes & \(0.0992 \pm 0.0245\) & \(0.1177 \pm 0.0305\) \\
\hline
Node-based CVT & \(0.1246 \pm 0.0244\) & \(0.1472 \pm 0.0302\) \\
\hline
Edge-based CVT & \(0.1901 \pm 0.0358\) & \(0.2367 \pm 0.0495\) \\
\hline
\end{tabular}
\end{center}
\end{table}

\subsection{Low-Level Path Planning}
We use an interior-point optimization algorithm, IPOPT~\cite{ipopt}, to solve  \Cref{eq:optimization}. Since IPOPT is a gradient-based optimizer, we provide the gradient of the objective function and the Jacobians of the constraints. These derivatives are computed using the JAX automatic differentiation library~\cite{jax2018github}.

\begin{figure*}[!htbp]
  \centering
  \includegraphics[scale=0.6, trim=0 0 0 30, clip]{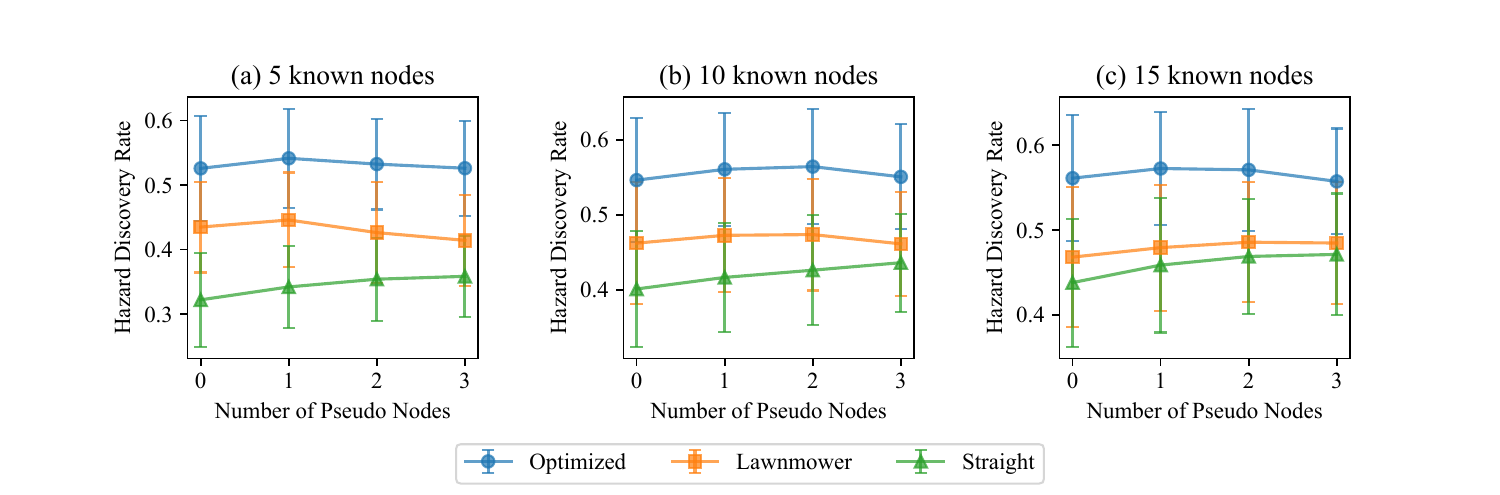}
  \caption{The average hazard discovery rate across 100 random trials is shown for three different path planning methods: optimized (blue circles), lawnmower (orange squares), and straight-line (green triangles). The results are presented for varying numbers of pseudo-nodes. Error bars represent one standard deviation across the trials. As shown, the optimized path consistently outperforms the other methods.}
  \label{fig:num_hazards_found}
\end{figure*}

To initialize the optimizer, we fit a B-spline to a ``lawnmower''-style path whose total length matches the path budget. This path is constructed within a rectangular region whose long axis lies along the straight-line segment connecting the two nodes, and whose short axis is perpendicular to that segment. The rectangle is centered along the segment between the two nodes, and its width is chosen so that the total length of the resulting back-and-forth path matches the path budget. The ``rungs'' of the pattern sweep across the rectangle, oriented perpendicular to the line between the nodes. The spacing between rungs is chosen based on the agent’s sensing range; specifically, the distance at which the probability of detection drops to $10\%$. An example of these paths can be seen in Figure~\ref{fig:planned_paths}(b). For comparison, we also evaluate unoptimized lawnmower paths of length $B_k$ as well as straight-line paths connecting the nodes. \Cref{fig:planned_paths} shows an example of each path type. As illustrated, the optimized paths tend to allocate more time to regions with a higher probability of containing hazards.

We ran the same 100 random trials described in the previous subsection, placing 5, 10, and 15 known hazards in the region. Additionally, we placed 50 unknown hazards, according to the distribution defined in \Cref{eq:combined_prior}, using $P_h = 0.3$ and $\lambda = 0.00015$. We selected 0, 1, 2, or 3 pseudo-nodes for each trial and planned routes accordingly. For each route, we generated paths using our optimization algorithm, as well as the lawnmower and straight-line baseline methods. We then allow the agent to traverse the path taking measurements according to the model in Equation~\eqref{eq:measurement_model} at a frequency of $1/\Delta_s = 0.1$ seconds, and $\beta=0.002$. The average number of unknown hazards discovered by the agent following each path is shown in \Cref{fig:num_hazards_found}.

As seen in the Figure~\ref{fig:num_hazards_found}, the optimized paths consistently outperform both the lawnmower and straight-line paths across all cases. The straight-line paths typically underutilize the path budget and thus offer limited exploration capability. The naive lawnmower paths do not account for the hazard prior distribution, resulting in the allocation of time inefficiently across the region. However, while pseudo-nodes contribute to improved coverage, the optimized results show limited gains in discovering unknown hazards. This is partly due to the static simulation setup, where no replanning occurs, and new hazards do not appear during execution. Addressing these limitations—such as enabling dynamic replanning—could further enhance the effectiveness of pseudo-nodes and is worth exploring in future work.


\section{Conclusion and Future Work}
This study presents a novel bi-level framework considering hazard exploration in UAV-based monitoring. Our approach efficiently balances exploration and exploitation in dynamic environments by integrating high-level route optimization with low-level path planning. The introduction of an edge-based CVT for pseudo-node placement significantly improves coverage, while Bayesian inference-based path planning maximizes information gain within constrained budgets. Simulation results validate the effectiveness of our method in improving hazard monitoring in uncertain environments.

For future work, we aim to enhance the path-planning phase by dynamically adjusting the path budget in response to real-time exploration feedback, enabling more effective replanning. We also plan to develop an iterative framework that integrates routing and path planning for continuous-monitoring scenarios, where both are updated based on new hazard data to improve responsiveness to evolving conditions.





\bibliographystyle{IEEEtran}
\bibliography{reference}  

\end{document}